\documentclass[reqno,11pt]{amsart} 
\usepackage{amsmath}
\usepackage{amssymb}
\usepackage{amsthm}
\usepackage{verbatim}
\usepackage{xcolor}
\usepackage{graphics}

\newcommand\NoBlackBoxes{\global\overfullrule0pt}

\NoBlackBoxes
\theoremstyle{plain}

\begin{document}

\title{Richter's local limit theorem, \\ its refinement,
and related results}

\author{S. G. Bobkov$^{1,4}$}
\thanks{1) 
School of Mathematics, University of Minnesota, Minneapolis, MN, USA,
bobkov@math.umn.edu. 
}

\author{G. P. Chistyakov$^{2,4}$}
\thanks{2) Faculty of Mathematics,
Bielefeld University, Germany,
chistyak@math-uni.bielefeld.de.
}

\author{F. G\"otze$^{2,4}$}
\thanks{2) Faculty of Mathematics,
Bielefeld University, Germany,
goetze@math-uni.bielefeld.de.
}

\thanks{4) Research supported by the NSF grant DMS-2154001 and 
the GRF – SFB 1283/2 2021 – 317210226 }

\subjclass[2010]
{Primary 60E, 60F} 
\keywords{Central limit theorem, local limit theorem, characteristic functions} 
\dedicatory{\centerline{In memoriam Gennadiy P. Chistyakov  \; 
*May 1, 1945 \;\;  $\dagger$December 30, 2022} }
	
\begin{abstract}
We give a detailed exposition of the proof of Richter's local limit 
theorem in a refined form, and establish the stability of the remainder
term in this theorem under small perturbations of the underlying distribution
(including smoothing). We also discuss related quantitative bounds 
for characteristic functions and Laplace transforms.
\end{abstract}

\maketitle
\markboth{S. G. Bobkov, G. P. Chistyakov and F. G\"otze}{Richter's local limit theorem}

\def\theequation{\thesection.\arabic{equation}}
\def\E{{\mathbb E}}
\def\R{{\mathbb R}}
\def\C{{\mathbb C}}
\def\P{{\mathbb P}}
\def\Z{{\mathbb Z}}
\def\S{{\mathbb S}}
\def\I{{\mathbb I}}
\def\T{{\mathbb T}}

\def\s{{\mathbb s}}

\def\G{\Gamma}

\def\Ent{{\rm Ent}}
\def\var{{\rm Var}}
\def\Var{{\rm Var}}
\def\V{{\rm V}}

\def\H{{\rm H}}
\def\Im{{\rm Im}}
\def\Tr{{\rm Tr}}
\def\s{{\mathfrak s}}

\def\k{{\kappa}}
\def\M{{\cal M}}
\def\Var{{\rm Var}}
\def\Ent{{\rm Ent}}
\def\O{{\rm Osc}_\mu}

\def\ep{\varepsilon}
\def\phi{\varphi}
\def\vp{\varphi}
\def\F{{\cal F}}

\def\be{\begin{equation}}
\def\en{\end{equation}}
\def\bee{\begin{eqnarray*}}
\def\ene{\end{eqnarray*}}

\thispagestyle{empty}

\vskip5mm
\section{{\bf Introduction and Formulation of the Results}}
\setcounter{equation}{0}

\vskip2mm
\noindent
Let $(X_n)_{n \geq 1}$ be independent copies of a random variable
$X$ with mean $\E X = 0$ and variance $\Var(X) = 1$.
Throughout, we assume without mentioning that the normalized sum
$$
Z_n = \frac{X_1 + \dots + X_n}{\sqrt{n}}
$$
has a bounded density $p_{n_0}$ for some $n = n_0$, that is, 
$p_{n_0}(x) \leq M$ for all $x \in \R$ with some constant $M$.
Then all $Z_n$ with $n \geq 2n_0$ have continuous bounded 
densities $p_n(x)$. An asymptotic behavior of these densities 
describing their closeness to the normal density
$$
\varphi(x) = \frac{1}{\sqrt{2\pi}}\,e^{-x^2/2}, \quad x \in \R,
$$
is governed by several local limit theorems. First of all, there is
a uniform local limit theorem due to Gnedenko
$$
\sup_x |p_n(x) - \varphi(x)| \rightarrow 0 \quad {\rm as} \
n \rightarrow \infty.
$$
Under higher order moment assumptions, say if $\E\,|X|^m < \infty$ 
for an integer $m \geq 3$, this statement may be considerably 
sharpened in the form of a non-uniform local limit theorem
\be
\sup_x \, (1 + |x|^m)\,|p_n(x) - \varphi_m(x)| =
o\big(n^{-\frac{m-2}{2}}\big),
\en
where $\varphi_m$ denotes the Edgeworth correction of $\varphi$
of order $m$ (cf. \cite{G-K}, \cite{Pe3}, \cite{Pe4}). In various 
applications, this relation is typically effective in the range 
$|x| \leq \sqrt{c\log n}$, since then the ratio $p_n(x)/\varphi(x)$ 
remains close to 1 (for a suitable $c$). For example, (1.1) is 
crucial in the study of rates in the entropic central limit theorem, 
rates for R\'enyi divergences of finite orders and for 
the relative Fisher information (\cite{B-C-G2}-\cite{B-C-G4}).

As for larger regions, the asymptotic behavior of
$p_n(x)$ is governed by the following remarkable theorem due 
to Richter \cite{R}, assuming the finiteness
of an exponential moment for the random variable $X$.

\vskip5mm
{\bf Theorem 1.1.} {\sl Suppose that, for some $b>0$,
\be
\E\,e^{b |X|} < \infty.
\en
Then, for $x = o(\sqrt{n})$, the densities of $Z_n$ admit 
the representation
\be
\frac{p_n(x)}{\varphi(x)} = \exp\Big\{\frac{x^3}{\sqrt{n}} \, 
\lambda\Big(\frac{x}{\sqrt{n}}\Big)\Big\} \,
\Big(1 + O\Big(\frac{1 + |x|}{\sqrt{n}}\Big)\Big).
\en
Here, $\lambda(\tau)$ represents
an analytic function in some neighborhood of zero.
}

\vskip5mm
It was shown by Amosova \cite{A} that the condition (1.2)
is necessary for the existence of a representation like (1.3) 
in the region $|x| = o(\sqrt{n})$ with some analytic 
function $\lambda$.

The function $\lambda$ in (1.3) is representable as a power series,
called the Cram\'er series,
\be
\lambda(\tau) = \sum_{k=0}^\infty \lambda_k \tau^k,
\en
which is absolutely convergent in some disc $|\tau| < \tau_0$ of 
the complex plane. It has appeared in the work by Cram\'er \cite{C} 
in a similar representation for the ratio of the tails of distribution 
functions of $Z_n$ and the standard normal law 
(cf. also \cite{F}, \cite{Pe1}, \cite{Pe2}).

Let us also mention that (1.3) is stated in Richter's 
work in a slightly different form with $O(\frac{|x|}{\sqrt{n}})$ 
in the last brackets and for $|x|>1$. A similar result is proved 
in the book by Ibragimov and Linnik \cite{I-L} under 
the assumption that $X$ has a bounded continuous density.

As a consequence of (1.3), one immediately obtains, for example,
that
\be
\frac{p_n(x)}{\varphi(x)} \rightarrow 1 \quad {\rm as} \ 
n \rightarrow \infty
\en
uniformly in the region $|x| = o(n^{1/6})$. In the region
$c_0 n^{1/6} \leq |x| \leq c_1 n^{1/2}$, the behavior 
may be quite different, 
and in order to describe it, the appearance of the term
$O(\frac{1+|x|}{\sqrt{n}})$ in (1.3) is non-desirable.
The purpose of this paper is to give a detailed exposition
of the proof of Theorem 1.1, clarifying the meaning of
the leading coefficient in (1.4) and replacing this term
with an $n$-depending quantity. We basically employ the
tools of \cite{I-L} and derive the following refinement.

\vskip5mm
{\bf Theorem 1.2.} {\sl Let the conditions of Theorem $1.1$
be fulfilled, and $n \geq 2n_0$.
There is a constant $\tau_0 > 0$ with the following property.
With $\tau = x/\sqrt{n}$, we have for $|\tau| \leq \tau_0$
\be
\frac{p_n(x)}{\varphi(x)} = e^{n\tau^3 \lambda(\tau) - \mu(\tau)}\,
\big(1 + O(n^{-1} (\log n)^3)\big),
\en
where $\mu(\tau)$ is an analytic function in $|\tau| \leq \tau_0$
such that $\mu(0)=0$.
}

\vskip5mm
As we will see in Section 5,
\bee
\lambda(\tau) 
 & = &
\frac{1}{m!}\,\gamma_m \tau^{m-3} + O(|\tau|^{m-2}), \\
\mu(\tau) 
 & = &
\frac{1}{2(m-2)!}\,\gamma_m \tau^{m-2} + O(|\tau|^{m-1}),
\ene
where $\gamma_m$ ($m \geq 3$) is the first non-zero cumulant 
of the random variable $X$ (assuming that it is not normal).
Equivalently, $m$ is the smallest positive integer such that
$\E X^m \neq \E Z^m$, where $Z$ is a standard normal 
random variable, in which case
$$
\gamma_m = \E X^m - \E Z^m.
$$

With this refinement, it should be clear that the relation (1.5) 
holds true uniformly over all $x$ in the potentially larger region
$$
|x| \leq \ep_n\,n^{\frac{1}{2} - \frac{1}{m}} \quad
(\ep_n \rightarrow 0).
$$ 
For example, if the distribution of $X$ is symmetric about the origin,
then $\gamma_3 = 0$, so that necessarily $m \geq 4$.

Another consequence of (1.6), which cannot be obtained on the
basis of (1.3), is needed in the study of the central limit theorem (CLT) 
with respect to the R\'enyi divergence of infinite order (including 
the rate of convergence). Let us recall that the R\'enyi divergence 
of a finite order $\kappa > 0$ from the distribution of $Z_n$ to
the standard normal law is defined by
$$
D_\kappa(p_n||\varphi) = \frac{1}{\kappa - 1} \log 
\int_{-\infty}^\infty \Big(\frac{p_n(x)}{\varphi(x)}\Big)^\kappa\,
\varphi(x)\, dx.
$$
As a function of $\kappa$, it is non-decreasing, representing a strong 
distance-like quantity. In the range $0<\kappa<1$, it is metrically 
equivalent to the total variation, that is, $L^1$-distance between 
$p_n$ and $\varphi$. The case $\kappa=1$ corresponds to the relative entropy
(Kullback-Leibler's distance)
$$
D(p_n||\varphi) = \lim_{\kappa \rightarrow 1} D_\kappa(p_n||\varphi) = 
\int_{-\infty}^\infty p_n(x)\,\log \frac{p_n(x)}{\varphi(x)}\,dx,
$$
and another important case
$\kappa = 2$ leads to the function of the $\chi^2$-Pearson distance.
So far, information-theoretic CLT's of the form
$D_\kappa(p_n||\varphi) \rightarrow 0$ as $n \rightarrow \infty$
have been completely characterized in terms of the distribution
of $X$ (that is, in the i.i.d. situation). However, such a statement
remains fully open for the limit distance
$$
D_\infty(p_n||\varphi) = \lim_{\kappa \rightarrow \infty} 
D_\kappa(p_n||\varphi) =
\sup_x\, \log \, \frac{p_n(x)}{\varphi(x)}.
$$
Equivalently, the problem is to find conditions under which
the related quantity (the Tsallis distance of infinite order) 
$$
T_\infty(p_n||\varphi) = 
\sup_x\, \frac{p_n(x) - \varphi(x)}{\varphi(x)}
$$
tends to zero for growing $n$ (note that one may not put the absolute 
values sign, since all $p_n$ may be compactly supported).

As a first natural step towards this variant of the CLT, we
consider the problem of the convergence for the restricted
Tsallis distance with the above suprema taken over
growing intervals $|x| = O(\sqrt{n})$. With this in mind, 
Theorem 1.2 allows one to prove the following:

\vskip5mm
{\bf Corollary 1.3.} {\sl Under the conditions of Theorem $1.1$, 
suppose that $m$ is even, $m \geq 4$, and $\gamma_m < 0$. There exist 
constants $\tau_0 > 0$ and $c>0$ with the following property.
If $|\tau| \leq \tau_0$, $\tau = x/\sqrt{n}$, then
\be
\frac{p_n(x) - \varphi(x)}{\varphi(x)} \leq \frac{c(\log n)^3}{n}.
\en
}

\vskip2mm
Here, the condition about cumulants is fulfilled, for example, when
the random variable $X$ is strongly subgaussian in the sense that
\be
\E\,e^{tX} \leq e^{t^2/2} \quad {\rm for \ all} \ t \in \R
\en
(recall that that $\E X^2 = 1$, while the condition $\E X = 0$ is necessary). 
This interesting class of probability distributions is rather rich, and 
we refer the reader to \cite{B-C-G5} for discussions and various examples. 
Our main motivation stemmed from the fact that the strong subgaussianity 
is necessary for the convergence $D_\infty(p_n||\varphi) \rightarrow 0$.
As we have recently learned, (1.8) had previously appeared under the name
``sharp subgaussianity" in the work by Guionnet and Husson \cite{G-H} 
for a completely different reason as a condition to have LDPs for the
largest eigenvalue of Wigner matrices with the same rate function as 
in the case of Gaussian entries.

One important issue which is not addressed in the formulation
of Theorem 1.2 is how one can control the involved constant
in the $O$-remainder term in (1.6). To better quantify this
asymptotic representation, we actually prove the following statement
using the same analytic functions $\lambda(\tau)$ and $\mu(\tau)$.

\vskip5mm
{\bf Theorem 1.4.} {\sl Assume that $\E\,e^{\alpha |X|} \leq 2$
$(\alpha>0)$. There exist absolute 
positive constants $C$ and $c$ such that, whenever $n \geq n_1$
and $\tau = x/\sqrt{n}$, $|\tau| \leq \tau_0$, we have
\be
\frac{p_n(x)}{\varphi(x)} = e^{n\tau^3 \lambda(\tau) - \mu(\tau)}
\Big(1 + B\alpha^{-6}\, \frac{(\log n)^3}{n}\Big),
\en
where $|B| \leq C$ and
$$
n_1 = CM^4 n_0^2\, \alpha^{-12}, \quad 
\tau_0 = \frac{c\alpha^3}{M^2 n_0}.
$$
}

\vskip2mm
This statement should be useful in applications
to smoothed distributions in order to guarantee that the
constant in the remainder term may be chosen to be common for all
distributions under consideration.

In order to make the proofs/arguments more transparent and self-contained,
we include a short review of various related results -- partly
technical, but often interesting in themselves -- about maxima of 
densities, analytic characteristic functions and log-Laplace transforms. 
The rest of the paper is organized as follows. In Sections 2 we recall
basic properties of the maximum of convolved densities and then
develop their applications to bounding restricted integrals of powers of 
characteristic functions (Section 3). In section 4, we discuss 
behavior of analytic characteristic functions near the origin.
Section 5 is devoted to the so-called saddle points and associated 
Taylor expansions for the log-Laplace transforms. Here we also analyze
the functions $\lambda(\tau)$ and $\mu(\tau)$.
Sections 6 deals with contour integration needed to establish preliminary
representations for $p_n(x)$. Final steps in the proof of Theorem 1.4
are made in Section 7. The proof of Corollary 1.3 is postponed to Section 8.

\vskip5mm
\section{{\bf Maximum of Convolved Densities}}
\setcounter{equation}{0}

\vskip2mm
\noindent
Convolved densities are known to have improved smoothing properties. 
First, let us emphasize the following general fact (which explains 
the condition $n \geq 2n_0$ mentioned before Theorem 1.1).

\vskip5mm
{\bf Proposition 2.1.} {\sl If independent random variables $\xi_1,\dots,\xi_m$
$(m \geq 2)$ have bounded densities, then the sum $S_m = \xi_1+ \dots +\xi_m$
has a bounded uniformly continuous density vanishing at infinity.
}

\vskip5mm
{\bf Proof.} Denote by $q_k$ the densities of $\xi_k$ and assume
that $q_k(x) \leq M_k$ for all $x \in \R$ with some constants
$M_k$ ($k \leq m$). By the Plancherel theorem, for the characteristic
functions $g_k(t) = \E\,e^{it\xi_k}$, we have
\bee
\int_{-\infty}^\infty |g_k(t)|^m\,dt
 & \leq &
\int_{-\infty}^\infty |g_k(t)|^2\,dt
 \, = \,
2\pi \int_{-\infty}^\infty q_k(x)^2\,dx \\
 & \leq &
2\pi \int_{-\infty}^\infty M_k q_k(x)\,dx \, = \, 2\pi M_k,
\ene
where we used the property $|g_k(t)| \leq 1$, $t \in \R$.
Hence, by H\"older's inequality, the characteristic function
$g(t) = g_1(t) \dots g_m(t)$ of $S_m$ is integrable and has
$L^1$-norm
\begin{eqnarray}
\int_{-\infty}^\infty |g(t)|\,dt
 & \leq & 
\Big(\int_{-\infty}^\infty |g_1(t)|^m\,dt\Big)^{1/m} \dots
\Big(\int_{-\infty}^\infty |g_m(t)|^m\,dt\Big)^{1/m} \nonumber \\
 & \leq &
2\pi \, (M_1 \dots M_m)^{1/m} \, < \, \infty.
\end{eqnarray}
One may conclude that the random variable $S_m$ 
has a bounded, uniformly continuous density expressed by 
the inversion Fourier formula
\be
q(x) = \frac{1}{2\pi} \int_{-\infty}^\infty e^{-itx} g(t)\,dt,
\quad x \in \R.
\en
Since $g$ is integrable, it also follows that $q(x) \rightarrow 0$
as $|x| \rightarrow \infty$ (by the Riemann-Lebesgue lemma).
\qed

\vskip2mm
Consider the functional
$$
M(\xi) = {\rm ess\,sup}_x \ q(x),
$$
where $\xi$ is a random variable with density $q$
(one may put $M(\xi) = \infty$ in all other cases). Since, by (2.2),
$$
q(x) \leq \frac{1}{2\pi} \int_{-\infty}^\infty |g(t)|\,dt
$$
for all $x \in \R$, the inequality (2.1) also implies that
\be
M(S_m) \leq (M(\xi_1) \dots M(\xi_m))^{1/m}.
\en
This shows in particular that $M(\xi)$ may not increase by adding 
to $\xi$ an independent random variable. However, the relation (2.3) 
does not correctly reflect  the behavior of $M(S_m)$ with respect 
to the growing parameter $m$, especially in the i.i.d. situation. 
A more precise statement is described in the following relation,
where the geometric mean of maxima is replaced with
the harmonic mean.

\vskip5mm
{\bf Proposition 2.2.} {\sl Given independent random variables
$\xi_k$, $1 \leq k \leq m$, one has
\be
\frac{1}{M(S_m)^2} \geq \frac{1}{2} \sum_{k=1}^m 
\frac{1}{M(\xi_k)^2}.
\en
}

\vskip2mm
This bound may be viewed as a counterpart of the entropy power 
inequality in Information Theory. It may be obtained by combining 
Rogozin's maximum-of-density theorem with Ball's bound on the
volume of slices of the cube. Namely, it 
was shown in \cite{Ro} that, if the values $M_k = M(\xi_k)$ are fixed,
$M(S_m)$ is maximized for $\xi_k$ uniformly distributed
in the intervals of length $1/M_k$. Of course, in this case $M(S_m)$
has a rather complicated structure as a function in variables 
$M_1,\dots,M_m$.

On the other hand, if
$T_m = a_1 \eta_1 + \dots + a_m \eta_m$, where
$\eta_k$ are independent and uniformly distributed in $(0,1)$,
and the coefficients satisfy $a_1^2 + \dots + a_m^2 = 1$,
then
\be
1 \leq M(T_m) \leq \sqrt{2},
\en
cf. \cite{Ball}. In geometric language, this is the same as saying
that $1 \leq |Q \cap H| \leq \sqrt{2}$, where $Q = (0,1)^m$ is the
unit cube, $H$ is an arbitrary hyperplane in $\R^m$ passing
through the center of the cube, and $|\,\cdot\,|$ stands for the
$(m-1)$-dimensional volume. To obtain (2.4), put 
$$
a_k = \frac{1}{aM_k}, \quad a^2 = \sum_{k=1}^m \frac{1}{M_k^2},
$$
so that, by the upper bound in (2.5),
\be
M\Big(\sum_{k=1}^m \frac{1}{M_k}\,\eta_k\Big) =
M\Big(a \sum_{k=1}^m a_k\,\eta_k\Big) = 
\frac{1}{a} M(T_m) \leq \frac{1}{a} \sqrt{2}.
\en
Since, by \cite{Ro}, $M(S_m)$ does not exceed the first term
in (2.6), we get $M(S_m) \leq \frac{1}{a} \sqrt{2}$, that is, (2.4). 

With this argument, this relation is mentioned in \cite{B-C}, 
where its multidimensional analog is derived by applying
the Hausdorff-Young inequality with best constants (due to
Beckner and Lieb). 

\vskip3mm
{\bf Remark 2.3.} Modulo a universal constant, the left inequality
in (2.5) may be extended to the more general setting. Namely, if 
a random variable $\xi$ has a density $q$ with a finite standard deviation 
$\sigma$, then
\be
M(\xi) \geq \frac{1}{12\,\sigma}.
\en
Here equality is attained for the uniform distribution on arbitrary
bounded intervals of the real line. This relation is well-known;  
as an early reference one can mention Statuljavičus \cite{St}, p.\,651,
where (2.7) is stated without proof. Since it is used below, 
let us include a short argument. For normalization, one may 
assume that $M(\xi) = 1$ and $\E \xi = 0$. In this
case, the tail function $H(x) = \P\{|\xi|\geq x\}$ has a
Lipschitz semi-norm at most 2, implying that
$H(x) \geq 1 - 2x$ for all $x \geq 0$. This gives
$$
\sigma^2 = \int_{-\infty}^\infty x^2\,q(x)\,dx =
2 \int_0^\infty x H(x)\,dx \geq 
2 \int_0^{1/2} x (1 - 2x)\,dx = \frac{1}{12}.
$$

\vskip5mm
\section{{\bf $L^p$-Norms of Characteristic Functions and Orlicz Norms}}
\setcounter{equation}{0}

\vskip2mm
\noindent
One useful consequence of (2.4) is
the next bound on $L^{2m}$-norms of characteristic functions.

\vskip5mm
{\bf Proposition 3.1.} {\sl If $g(t)$ is the characteristic function 
of a random variable $\xi$, then for any integer $m \geq 1$,
\be
\frac{1}{2\pi} \int_{-\infty}^\infty |g(t)|^{2m}\,dt \leq 
\frac{1}{\sqrt{m}}\, M(\xi).
\en
}

\vskip2mm
{\bf Proof.} We apply Proposition 2.2 to $2m$ summands 
$\xi_1,-\xi_1',\dots,\xi_m,-\xi_m'$, assuming that $\xi_k$, $\xi_k'$ are
independent copies of $\xi$. Introduce the symmetrized random variable 
$\tilde S_m = S_m - S_m'$, where $S_m'$ is an independent copy 
of $S_m$. By (2.4), we then get
$$
M(\tilde S_m) \leq \frac{1}{\sqrt{m}}\, M(\xi).
$$
In addition, $\tilde S_m$ has characteristic function $|g(t)|^{2m}$.
If $M(\xi)$ is finite, one may apply Proposition 2.1 and conclude 
that $\tilde S_m$ has a bounded continuous density $q_m(x)$ which
is vanishing at infinity. Moreover, $q_m(x)$ is maximized at $x = 0$, and 
its value at this point is described by the inversion formula (2.2) which gives
$$
M(\tilde S_m) = q_m(0) = 
\frac{1}{2\pi} \int_{-\infty}^\infty |g(t)|^{2m}\,dt.
$$
\qed

\vskip5mm
Using (2.3), one can obtain a similar relation, but without
the factor $1/\sqrt{m}$ in (3.1).

When $M(\xi)$ is finite and $m$ is large, this bound may be 
considerably sharpened asymptotically with respect to $m$
when restricting the integration to the regions $|t| \geq \ep > 0$.
Before making this precise, first let us note that,
since the random variable $\xi$ has a density, we have
\be
\delta_g(\ep) = \max_{|t| \geq \ep} |g(t)| < 1
\en
for all $\ep>0$. This holds by continuity of $g$, and since $|g(t)| < 1$ 
for all $t \neq 0$ (which is true for any non-lattice distribution), 
while $g(t)$ tends to zero as $t \rightarrow \infty$, by the 
Riemann-Lebesgue lemma. By the way, this property remains to hold 
in the more general situation, where the $m$-fold convolution of the 
distribution of $\xi$ with itself has a density (while the 
distribution of $\xi$ might be not absolutely continuous). 
Indeed, in that case, (3.2) may be applied to $g^m$, and it remains
to notice that this relation does not depend on $m$.

The property (3.2) may be quantified using, for example, the following 
observation due to Statuljavičus \cite{St}.

\vskip5mm
{\bf Proposition 3.2.} {\sl If a random variable $\xi$ has a bounded
density with $M = M(\xi)$ and finite variance $\sigma^2 = \Var(\xi)$, 
$\sigma>0$, then its characteristic function $g$ satisfies, for all $\ep > 0$,
\be
\delta_g(\ep) \leq 
\exp\Big\{-\frac{\ep^2}{96\,M^2\,(2\sigma \ep + \pi)^2} \Big\}.
\en
}

\vskip2mm
This relation may be extended to non-bounded
densities $q$, in which case the parameter $M$ should be replaced
with quantiles of the random variable $q(\xi)$. The moment condition
may also be removed, and instead it is sufficient to deal with
quantiles of $|\xi - \xi'|$, where $\xi'$ is an independent copy of $\xi$;
cf. \cite{B-C-G1} for details.

Returning to (3.1) and applying (3.3) with $\ep \leq 1$, we then have
$$
\int_{|t| \geq \ep} |g(t)|^{4m}\,dt 
 \, \leq \,
\delta_g(\ep)^{2m} \int_{-\infty}^\infty |g(t)|^{2m}\,dt
  \, \leq \,
\frac{2\pi M}{\sqrt{m}}\, \exp\Big\{-\frac{m\ep^2}{CM^2}\Big\}
$$
with some absolute constant $C$.
Thus, the resulting bound decays asymptotically fast in~$m$.

Let us derive a similar bound in the scheme of independent copies
$(X_n)_{n \geq 1}$ of the random variable $X$ with
$\Var(X) = 1$, assuming that the normalized sum $Z_n$
has a bounded density for $n = n_0$ with
$M = M(Z_{n_0})$. Consider the characteristic function
$f(t) = \E\,e^{itX}$.
We apply Propositions 3.1-3.2 with $\xi = X_1 + \dots + X_{n_0}$ 
in which case $g(t) = f(t)^{n_0}$ and $M(\xi) = \frac{1}{\sqrt{n_0}}\,M$.
Then, for any $1 \leq m \leq n/2n_0$, by (3.1),
\bee
\int_{|t| \geq \ep} |f(t)|^n\,dt 
 & = &
\int_{|t| \geq \ep} |f(t)|^{n - 2m n_0}\, |g(t)|^{2m}\,dt \\
 & \leq &
\delta_f(\ep)^{n-2mn_0} \int_{-\infty}^\infty |g(t)|^{2m}\,dt
 \, \leq \,
\frac{2\pi M}{\sqrt{m n_0}}\, \delta_f(\ep)^{n-2m n_0}.
\ene
If $n \geq 4n_0$, let us choose $m = [\frac{n}{4n_0}]$. Then 
$n - 2m n_0 \geq \frac{n}{2}$, while $m \geq \frac{n}{8 n_0}$,
and we arrive at
$$
\int_{|t| \geq \ep} |f(t)|^n\,dt \leq \frac{4\pi M}{\sqrt{2n}}\,
\delta_f(\ep)^{n/2}.
$$
By (3.3) with $\ep \leq 1$, we also have
$$
(\delta_f(\ep))^{n_0} \leq 
\exp\Big\{-\frac{\ep^2}{96\,M^2\,(2\sqrt{n_0} + \pi)^2} \Big\},
$$
which may be simplified to
$$
\delta_f(\ep) \leq 
\exp\Big\{-\frac{\ep^2}{96\,(2 + \pi)^2\, n_0 M^2} \Big\}.
$$
Combining the two bounds, one may summarize.

\vskip5mm
{\bf Corollary 3.3.} {\sl Let $\Var(X) = 1$, and suppose that $Z_n$ has
a density for $n = n_0$ bounded by $M$. Then, for all $0 < \ep \leq 1$ 
and $n \geq 4n_0$, the characteristic function $f$ of $X$ satisfies
\be
\int_{|t| \geq \ep} |f(t)|^n\,dt \leq \frac{4\pi M}{\sqrt{2n}}\,
\exp\Big\{-\frac{n\ep^2}{Cn_0 M^2}\Big\}, \quad C = 5200.
\en
}

\vskip5mm
\section{{\bf Behavior of Characteristic Functions near Zero}}
\setcounter{equation}{0}

\vskip2mm
\noindent
While the boundedness of the density is important to control
integrability properties of powers of the characteristic 
function of a random variable $X$, the condition (1.2) 
on the finiteness of an exponential moment of $X$ guarantees that
the characteristic function
$$
f(z) = \E\,e^{izX}, \quad z = t+iy, \quad t,y \in \R,
$$
is well-defined and analytic in the strip $|y| = |{\rm Re}(z)| < b$ 
of the complex plane. Equivalently, we will assume throughout
that, for some $\alpha > 0$,
\be
\E\,e^{\alpha |X|} \leq 2.
\en 
This parameter is more convenient to quantify the behavior of $f(z)$ near zero.

For example, using $x e^{-x} \leq e^{-1}$ ($x \geq 0$) and assuming that
$|y| \leq \frac{\alpha}{2}$, we then have
\bee
|f'(z)| 
 & = & 
|\,\E\,X e^{iz X}|\, \leq \,
\E\,|X|\,e^{|yX|} \, \leq \, \E\,|X|\,e^{\alpha |X|/2} \\
 & = &
\E\,|X|\,e^{-\alpha |X|/2} \,e^{\alpha|X|} \, \leq \, \frac{4}{\alpha e}.
\ene
Hence $|f(z) - 1| \leq \frac{4}{\alpha e}\, |z|$ (since $f(0) = 1$).
Thus, we obtain:

\vskip5mm
{\bf Lemma 4.1.} {\sl For all complex numbers $z$ in the disc
$|z| \leq \frac{\alpha}{2}$, 
$$
|f'(z)| \leq \frac{4}{\alpha e}, \quad |f(z) - 1| \leq \frac{2}{e}.
$$
}

\vskip2mm
This allows one to consider the log-Laplace transform
$$
K(z) = \log \E\,e^{zX} = \log f(-iz)
$$
as an analytic function in the disc $|z| \leq \frac{\alpha}{2}$.
Since it has derivative
$
K'(z) = -i\,\frac{f'(-iz)}{f(-iz)},
$
from Lemma 4.1 we get that in this disc
\be
|K'(z)| \leq \frac{6}{\alpha}, \quad 
|K(z)| \leq 3.
\en
One may also bound the derivatives of all orders.

\vskip5mm
{\bf Lemma 4.2.} {\sl For all complex numbers $z$ in the disc
$|z| \leq \frac{\alpha}{4}$, 
\be
|K^{(k)}(z)| \leq 3 k!\,\Big(\frac{4}{\alpha}\Big)^k, \quad
k = 1,2,\dots
\en
Moreover, if $\E X = 0$, $\E X^2 = 1$, then
\be
|K'''(z)| \leq \frac{8}{\alpha^3}, \quad |z| \leq \frac{\alpha}{16}.
\en
As a consequence,
\be
|K''(z) - 1| \leq \frac{1}{2}, \quad |z| \leq \frac{\alpha^3}{16}.
\en
}

\vskip2mm
Thus, these derivatives have at most a factorial growth
in absolute value with respect to the growing parameter $k$.
For the particular orders $k=2$ and $k=3$, and under our
moment assumptions, the bound (4.3) may be refined
in a smaller disc according to (4.4)-(4.5).

\vskip5mm
{\bf Proof.} To obtain (4.3), one may apply Cauchy's formula
$$
K^{(k)}(z) = \frac{k!}{2\pi i} \int_{|w-z| = r} 
\frac{K(w)}{w^{k+1}}\,dw
$$
with $r = \frac{\alpha}{4}$ together with the second bound in (4.2).

Turning to the refined bounds, note that in terms of the Laplace transform 
$L(z) = \E\,e^{zX}$, we have $K' = L'/L$ and
\be
K''' = \frac{L'''}{L} - 3\,\frac{L'' L'}{L^2} + 2\, \frac{L'^3}{L^3}.
\en

For $x\geq 0$ and $p = 1,2,3$, we use the elementary inequality 
$x^p e^{-x} \leq (p/e)^p$. Suppose that $|z| \leq (1 - c)\, \alpha$ 
with $\frac{1}{2} \leq c < 1$. Since $L^{(p)}(z) = \E\,X^p\, e^{zX}$, 
we then have
$$
|L^{(p)}(z)| \leq \E\,|X|^p\, e^{(1 - c)\,\alpha |X|}.
$$ 
Hence, by (4.1),
$$
|L^{(p)}(z)| \leq \,
\E\,|X|^p\,e^{-c \alpha |X|} \,e^{\alpha |X|} \, \leq \, 
2\,\Big(\frac{p}{c \alpha e}\Big)^p.
$$
In particular,
\be
|L'(z)| \leq \frac{2}{c\alpha e}, \quad 
|L''(z)| \leq \frac{8}{(c\alpha e)^2}, \quad 
|L'''(z)| \leq \frac{54}{(c\alpha e)^3},
\en
so
$$
|L(z) - 1| \leq \frac{2}{c\alpha e}\,|z| \leq 
\frac{2(1-c)}{c e}, \qquad
|L(z)| \geq 1 - \frac{2(1-c)}{c e}.
$$
Putting 
$q^{-1} = 1 - \frac{2(1-c)}{c e}$, from (4.6) it follows that
$$
|K'''(z)| \leq (c\alpha e)^{-3}\, (54\,q + 48\,q^2 + 16\,q^3).
$$
Choosing $c = 15/16$, we have
$q = (1 - \frac{2}{15\,e})^{-1} < 1.06$, and the last expression
becomes smaller than $8\alpha^{-3}$. Hence
$$
|K'''(z)| \leq \frac{8}{\alpha^3}, \quad |K''(z) - 1| \leq \frac{8}{\alpha^3}\,|z|
$$
for $|z| \leq \frac{\alpha}{16}$, where we used $K''(0)=1$.
The last inequality readily implies (4.5).
\qed

\vskip5mm
We shall now show that $|f(z)|$ is bounded away from 1
in a certain region near zero.

\vskip5mm
{\bf Lemma 4.3.} {\sl Let $\E X = 0$, $\E X^2 = 1$. For all
complex numbers $z = t+iy$ with 
$|t| \leq \frac{\alpha^3}{8}$, $|y| \leq \frac{1}{2}\,|t|$, we have
\be
|f(z)| \leq e^{-t^2/5}.
\en
}

\vskip2mm
{\bf Proof.} Using $f'(0) = 0$, $f''(0) = -1$,
one may start with an integral Taylor formula
$$
f(z) = 1 - \frac{1}{2}\,z^2 + 
\frac{1}{2}\,z^3 \int_0^1 f'''(s z) (1-s)^2\,ds.
$$
This equality is needed in the disc $|z| \leq r$ of radius 
$r = \frac{1}{8}\sqrt{\frac{5}{4}}\, \alpha^3$.
By the triangle inequality, we then have
\be
|f(z)| \leq 
\Big|1 - \frac{1}{2}\,z^2\Big| + \frac{A}{6}\,|z|^3, \quad
A = \max_{|z| \leq r} |f'''(z)|.
\en

First, let us check that
\be
\Big|1 - \frac{1}{2}\,z^2\Big| \leq 1 - \frac{1}{3}\,t^2,
\en
which actually holds in the larger region 
$|y| \leq \frac{1}{2}\,|t|$, $|t| \leq \frac{1}{4}$.
In this case, $t^2 - y^2 \geq \frac{3}{4}\,t^2$ and
$|ty| \leq \frac{1}{2}\,t^2$, implying
\bee
\Big|1 - \frac{1}{2}\,z^2\Big|^2
 & = &
1 - (t^2 - y^2) + \frac{1}{4}\,(t^2 - y^2)^2 + (ty)^2 \\
 & \leq &
1 - \frac{3}{4}\,t^2 + \frac{1}{2}\,t^4
 \, \leq \,
\Big(1 - \frac{1}{3}\,t^2\Big)^2,
\ene
where we used $|t| \leq \frac{1}{4}$ on the last step.
Thus, (4.10) follows.

Turning to the maximum in (4.9), one may apply the last bound
in (4.7) valid for $|z| \leq (1-c)\, \alpha$. Hence, we have
the constraint $(1-c) \alpha \geq r$, which is fulfilled
for the choice $c = 1 - \frac{1}{8}\sqrt{\frac{5}{4}}$
(due to $\alpha < 1$, cf. Lemma 4.4 below). In this case, we get
$$
|f'''(z)| \leq \frac{54}{(c\alpha e)^3} < \frac{4.3}{\alpha^3}.
$$
For $z = t+iy$, $|y| \leq \frac{1}{2}\,|t|$, we have 
$|z|^3 \leq (\frac{5}{4})^{3/2}\,|t|^3$,
and (4.9)-(4.10) therefore give
\bee
|f(z)|
 & \leq &
1 - \frac{1}{3}\,t^2 + 
\frac{4.3}{6 \alpha^3}\,\Big(\frac{5}{4}\Big)^{3/2}\,|t|^3 \\
 & \leq &
1 - \frac{1}{3}\,t^2 + 
\frac{4.3}{48}\,\Big(\frac{5}{4}\Big)^{3/2}\,t^2 \, \leq \,
1 - \frac{1}{5}\,t^2.
\ene
\qed

\vskip5mm
Finally, let us make a few remarks about the relationship between 
the conditions (1.2) and (4.1). When the random variable $X$ has 
a finite exponential moment, and $\alpha$ is optimal, then (4.1) 
becomes an equality. In this case, the quantity $\frac{1}{\alpha}$ 
represents the Orlicz norm of $X$ generated 
by the Young function $\psi(x) = e^{|x|} - 1$, $x \in \R$:
$$
\|X\|_{\psi} = 
\inf\big\{\lambda > 0: \E\,\psi(X/\lambda) \leq 1\big\}.
$$

If $\E X^2 = 1$, the parameter $\alpha$ may not be large, 
since the $L^2$-norm is dominated by the $L^\psi$-norm.
More precisely, using $x^2 e^{-x} \leq 4e^{-2}$ ($x \geq 0$), 
we have
$$
\alpha^2 = \E\, (\alpha X)^2 \leq 
4e^{-2}\, \E\,e^{\alpha |X|} = 8e^{-2},
$$ 
implying $\alpha \leq 2e^{-1} \sqrt{2} < 1.05$.
In fact, this bound may be sharpened.

\vskip5mm
{\bf Lemma 4.4.} {\sl If $\E X^2 = 1$, and $(4.1)$ holds, then 
$\alpha < 1$.
}

\vskip5mm
{\bf  Proof.} We may assume that $X \geq 0$ and then we 
need to show that $\E\,e^X > 2$. It is easy to check that
$x + \frac{1}{6}\,x^3 \geq ax^2$ for all $x \geq 0$ with
the optimal constant $a = \frac{2}{\sqrt{6}}$. Since 
$\E X^k \geq (\E X^2)^{k/2} = 1$ for $k \geq 2$, we get
\bee
\E\,e^X 
 & = &
1 + \frac{1}{2}\,\E X^2 + \E\,\Big(X + \frac{1}{6}\,X^3\Big) +
\sum_{k=4}^\infty \frac{1}{k!}\,\E X^k \\
 & \geq &
\frac{3}{2} + a + \sum_{k=4}^\infty \frac{1}{k!} \, = \,
e - \frac{7}{6} + \frac{2}{\sqrt{6}} \, > \, 2.36.
\ene
\qed

\vskip2mm
Note that if we start with a more general condition 
$B = \E e^{b |X|} < \infty$ as in Theorem 1.1, (4.1) is fulfilled 
for a certain constant $\alpha>0$. Indeed, if $B \leq 2$, then 
one may take $\alpha=b$. Otherwise,
$$
\E\,e^{\ep b|X|} \leq (\E\,e^{b|X|})^\ep \leq B^\ep = 2
$$
for $\ep = \frac{1}{\log_2(B)}$. Hence $\alpha = \ep b = \frac{b}{\log_2(B)}$
works as well. The two cases may be united by taking
\be
\alpha = \frac{b}{\log_2(\max(B,2))}.
\en

\vskip5mm
\section{{\bf Saddle Point and Taylor Expansions}}
\setcounter{equation}{0}

\vskip2mm
\noindent
Assume that $\E X = 0$, $\E X^2 = 1$, and $\E\,e^{\alpha |X|} \leq 2$
($\alpha>0$). Since the log-Laplace transform $K(z) = \log \E\,e^{zX}$ 
was defined as an analytic function in the disc $|z| \leq \frac{\alpha}{2}$ 
of the complex plane, it 
may be expanded as an absolutely convergent power series
$$
K(z) = \frac{1}{2}\,z^2 + 
\sum_{k=3}^\infty \frac{\gamma_k}{k!}\,z^k.
$$
Here, the coefficients $\gamma_k = K^{k)}(0)$ are called cumulants 
of $X$. Every $\gamma_k$ represents a certain polynomial in moments
of $X$ up to order $k$. In particular, $\gamma_3 = \E X^3$ and
$\gamma_4 = \E X^4 - 3$.

Similarly,
$$
K'(z) = z + \sum_{k=2}^\infty \frac{\gamma_{k+1}}{k!}\,z^k.
$$

The next object is important for contour integration.

\vskip5mm
{\bf Definition 5.1.} Given $\tau \in \C$, a saddle point is a solution
$z_0 = z_0(\tau)$ of the equation
\be
K'(z) = \tau.
\en

\vskip5mm
Thus, a saddle point is the solution of
\be
\tau = z  +
\sum_{k=2}^\infty \frac{\gamma_{k+1}}{k!}\,z^k.
\en

\vskip5mm
{\bf Proposition 5.2.} {\sl In the disc $|\tau| \leq \frac{\alpha^3}{32}$,
the equation $(5.1)$ has a unique solution $z_0(\tau)$. Moreover, it 
represents an injective analytic function satisfying $z_0'(0)=1$ and
\be
|z_0(\tau)| \leq 2 \tau \leq \frac{\alpha^3}{16}, \quad
|\tau| \leq \frac{\alpha^3}{32}.
\en
}

\vskip2mm
{\bf Proof.}
Let us use (5.2) as the definition of the analytic 
function $\tau = K'(z)$. If $\tau$ is sufficiently small, 
say $|\tau| \leq \tau_0$, this equality
may be inverted as a power series in $\tau$, 
$$
z = z_0(\tau) = \tau - \frac{\gamma_3}{2}\,\tau^2 + 
\frac{3\gamma_3^2 - \gamma_4}{6}\,\tau^3 + \dots
$$
Let us indicate an explicit expression for $\tau_0$ in the form
of a positive function of $\alpha$. 

By Lemma 4.2, cf. (4.4),
\be
|\tau'(z) - 1| \leq \frac{8}{\alpha^3}\,|z|, \quad 
|z| \leq \frac{\alpha}{16}.
\en
One may use this relation for $|z| \leq \frac{\alpha^3}{16}$,
since $\alpha<1$, by Lemma 4.4. Given two points $z_1$ and $z_2$ 
in the disc $|z| \leq \frac{\alpha^3}{16}$, define the path 
$z_t = (1-t) z_1 + t z_2$ connecting these points. We have
\be
\tau(z_2) - \tau(z_1) =
(z_2 - z_1) \Big(1 + \int_0^1 (\tau'(z_t) - 1)\,dt\Big),
\en
implying
$$
|\tau(z_2) - \tau(z_1)| \geq 
|z_2 - z_1| \Big(1 - \int_0^1 |\tau'(z_t) - 1|\,dt\Big).
$$
Since $|z_t| \leq \frac{\alpha^3}{16}$, it follows from (5.4) that
$$
|\tau(z_2) - \tau(z_1)| \geq \frac{1}{2}\,|z_2 - z_1|.
$$

As a consequence, the map $z \rightarrow \tau(z)$ is injective
in the disc $|z| \leq \frac{\alpha^3}{16}$. In addition, since 
$\tau(0)=0$, we have
\be
|\tau(z)| \geq \frac{1}{2}\,|z|.
\en
Therefore, the image of the circle $|z| = \frac{\alpha^3}{16}$
under this map represents a closed curve on the complex plane 
outside the circle 
$|\tau| = \frac{\alpha^3}{32}$. Since the image of the disc 
$|z| \leq \frac{\alpha^3}{16}$ under $\tau$ is a connected set, 
while $\tau(0)=0$, this set must contain the disc 
$|\tau| \leq \frac{\alpha^3}{32}$. 
Thus, the inverse map $z_0(\tau) = \tau^{-1}$ is well-defined and 
represents a holomorphic injective function in
$|\tau| \leq \frac{\alpha^3}{32}$ satisfying (5.3), by (5.6), and
$z_0'(0)=1$, by (5.4). Hence, one may take 
$\tau_0 = \frac{\alpha^3}{32}$. 

In addition, $z_0(\tau)$ takes real values for real $\tau$.
Indeed, since all cumulants are real numbers, $\tau(z)$
is real for real $z$, so is the inverse function $z_0$. Also, by (5.5), 
$$
\tau = z_0(\tau) \Big(1 + \int_0^1 (\tau'(t z_0(\tau)) - 1)\,dt\Big),
$$
which shows that $z_0(\tau)>0$ as long as 
$0 < \tau \leq \frac{\alpha^3}{32}$ (since the expression under 
the integral sign is a real-valued function whose absolute value
does not exceed $1/2$).
\qed

\vskip5mm
It is natural to determine the leading term in the Taylor expansion 
for $z_0(\tau)$ when expanding this function as a power series
in $\tau$. Assuming that $X$ is not normal, let $\gamma_m$ ($m \geq 3$) 
be the first non-zero cumulant of $X$. Then, as $|z| \rightarrow 0$,
$$
K(z) = \frac{1}{2}\,z^2 + \frac{\gamma_m}{m!}\,z^m +
O(|z|^{m+1}),
$$
so that
\be
K'(z) = z + \frac{\gamma_m}{(m-1)!}\,z^{m-1} +
O(|z|^m)
\en
and
\be
K''(z) = 1 + \frac{\gamma_m}{(m-2)!}\,z^{m-2} +
O(|z|^{m-1}).
\en
Since, $z_0(\tau) = \tau + O(|\tau|^2)$ as $\tau \rightarrow 0$, 
cf. Proposition 5.2, we get from (5.7)
\bee
\tau 
 & = &
K'(z_0(\tau)) \, = \,
z_0(\tau) + \frac{\gamma_m}{(m-1)!}\,z_0(\tau)^{m-1} +
O(|z_0(\tau)|^m) \\
 & = &
z_0(\tau) + \frac{\gamma_m}{(m-1)!}\,\tau^{m-1} +
O(|\tau|^m).
\ene
Therefore,
\be
z_0(\tau) = \tau - \frac{\gamma_m}{(m-1)!}\,\tau^{m-1} +
O(|\tau|^m).
\en

Next, let us write down the Taylor expansion around the point 
$z_0 = z_0(\tau)$:
\be
K(z) - \tau z = K(z_0) - \tau z_0 + \sum_{k=2}^\infty
\frac{\rho_k}{k!}\,(z-z_0)^k, \quad
\rho_k = K^{(k)}(z_0).
\en
Here, we used the property that the function $K(z) - \tau z$
has derivative $K'(z_0) - \tau = 0$ at the saddle point $z=z_0$.
Thus, the linear term in (5.10) corresponding to $k=1$ is vanishing.
As for the free term corresponding to $k=0$, we have
$$
K(z_0) = \sum_{k=2}^\infty \frac{\gamma_k}{k!}\,z_0^k
$$
and
$$
\tau z_0 = z_0 K'(z_0) = 
\sum_{k=2}^\infty \frac{\gamma_k}{(k-1)!}\,z_0^k.
$$
Hence
\be
K(z_0) - \tau z_0 = - \sum_{k=2}^\infty 
\frac{k-1}{k!}\,\gamma_k z_0^k =
-\frac{1}{2}\,z_0^2 - \frac{1}{3}\,\gamma_3 z_0^3 + \dots
\en
Using (5.9) and (5.11), we actually have
\bee
K(z_0) - \tau z_0 
 & = &
-\frac{1}{2}\,z_0^2 - \frac{m-1}{m!}\,\gamma_m z_0^m + \dots \\
  & = &
-\frac{1}{2}\,\Big(\tau - \frac{\gamma_m}{(m-1)!}\,\tau^{m-1} +
O(|\tau|^m)\Big)^2 \\
 & & - \
\frac{m-1}{m!}\,\gamma_m 
\Big(\tau - \frac{\gamma_m}{(m-1)!}\,\tau^{m-1} +
O(|\tau|^m)\Big)^m + \dots
\ene
which is simplified to
\begin{eqnarray}
K(z_0) - \tau z_0
 & = &
-\frac{1}{2}\,\tau^2 + \frac{1}{m!}\,\gamma_m \tau^m + 
O(|\tau|^{m+1}) \nonumber \\
 & = &
-\frac{1}{2}\,\tau^2 + \tau^3 \lambda(\tau).
\end{eqnarray}
Thus, applying Proposition 5.2 and recalling that
$K(z)$ is analytic in $|z| \leq \frac{\alpha}{2}$ (Lemma 4.1), we obtain:

\vskip5mm
{\bf Proposition 5.3.} {\sl The function
$$
\lambda(\tau) = 
\frac{1}{\tau^3}\,\Big(K(z_0(\tau)) - \tau z_0(\tau) + \frac{1}{2}\,\tau^2\Big)
$$
is well-defined and analytic in the disc $|\tau| \leq \frac{\alpha^3}{32}$.
Moreover, as $\tau \rightarrow 0$,
\be
\lambda(\tau) = \frac{1}{m!}\,\gamma_m \tau^{m-3} + 
O(|\tau|^{m-2}).
\en
}

\vskip2mm
{\bf Definition 5.4.} Being an analytic function, $\lambda(\tau)$ is 
represented as a power series in the disc 
$|\tau| \leq \frac{\alpha^3}{32}$. It is called Cram\'er's series.

\vskip5mm
It follows that $\lambda(\tau)$ is bounded for small $\tau$, but we will
need to quantify this property in terms of the parameter $\alpha$.
Recall that $\E X = 0$, $\E X^2 = 1$ and $\E\,e^{\alpha |X|} \leq 2$ ($\alpha>0$).

\vskip5mm
{\bf Proposition 5.5.} {\sl We have
\be
|\lambda(\tau)| \leq 700\,\alpha^{-3}, \quad |\tau| \leq \frac{\alpha^3}{64}.
\en
}

\vskip2mm
{\bf Proof.} Proposition 5.2 allows us to apply Cauchy's formula, which yields
$$
z_0''(\tau) = \frac{4}{2\pi i} 
\int_{|\xi - \tau| = r} \frac{z_0(\xi)}{(\xi - \tau)^3}\,d\xi
$$
with $r = \frac{\alpha^3}{64}$. Moreover, the latter implies, by (5.3),
\be
|z_0''(\tau)| \, \leq \, \frac{4}{r^2}\, \max_{|\xi - \tau| = r} |z_0(\xi)|
 \, \leq \, \frac{4}{r^2}\cdot \frac{\alpha^3}{16} 
 \, = \, 2^{12}\,\alpha^{-3}.
\en

Next, we note that, by Definition 5.1 of the saddle point, cf. (5.1),  the function 
$$
\psi(\tau) = K(z_0(\tau)) - \tau z_0(\tau) + \frac{1}{2}\,\tau^2
$$
has the first three derivatives
\bee
\psi'(\tau) 
 & = &
K'(z_0(\tau)) z_0'(\tau) - z_0(\tau) - \tau z_0'(\tau) + \tau \, = \,
\tau - z_0(\tau), \\
\psi''(\tau) 
 & = &
1 - z_0'(\tau), \\
\psi'''(\tau) 
 & = & 
- z_0''(\tau).
\ene
Since $\psi(0) = \psi'(0) = \psi''(0) = 0$, we may apply 
the Taylor integral formula together with (5.15) to conclude that
\bee
|\psi(\tau)|
 & \leq &
\frac{|\tau|^3}{6}\,\max_{|\xi| \leq |\tau|} |\psi'''(\xi)| \\
 & \leq &
\frac{|\tau|^3}{6}\,\max_{|\xi| \leq |\tau|} |z_0''(\xi)| \, \leq \,
\frac{1}{6}\cdot 2^{12}\,\alpha^{-3}.
\ene
As $\psi(\tau) = \tau^3 \lambda(\tau)$, the relation (5.14) follows.
\qed

\vskip5mm
Let us now introduce another analytic function which appears
in the representation (1.6) of Theorem 1.2.

\vskip5mm
{\bf Proposition 5.6.} {\sl The function
$$
\mu(\tau) = \frac{1}{2}\,\log K''(z_0(\tau))
$$
is well-defined and analytic in the disc $|\tau| \leq \frac{\alpha^3}{32}$.
Moreover, as $\tau \rightarrow 0$,
\be
\mu(\tau) = \frac{1}{2(m-2)!}\,\gamma_m \tau^{m-2} + 
O(|\tau|^{m-1}).
\en
}

\vskip2mm
{\bf Proof.} By (5.3), $|z_0(\tau)| \leq \frac{\alpha^3}{16}$.
Hence, by Lemma 4.2, $K''(z_0(\tau))$ takes values in the disc
with center at 1 of radius $1/2$. Thus, the principal value of
$\log K''(z_0(\tau))$ is well-defined and represents
an analytic function in $|\tau| \leq \frac{\alpha^3}{32}$.
Moreover, by (5.8)-(5.9),
\bee
K''(z_0(\tau)) 
 & = &
1 + \frac{\gamma_m}{(m-2)!}\,z_0(\tau)^{m-2} +
O(|z_0(\tau)|^{m-1}) \\
 & = &
1 + \frac{\gamma_m}{(m-2)!}\,\tau^{m-2} +
O(|\tau|^{m-1}).
\ene
Taking the logarithm of this expression, we arrive at (5.16).
\qed

\vskip5mm
Let us also mention that the function $K(z)$ is convex
and has a positive second derivative on the real line, more
precisely -- on the interval where it is finite. Hence $\mu(\tau)$
is real-valued for real $\tau$.

\vskip5mm
\section{{\bf Contour Integration}}
\setcounter{equation}{0}

\vskip2mm
\noindent
Let $(X_n)_{n \geq 1}$ be independent copies of a random variable
$X$ with $\E X = 0$, $\Var(X) = 1$, and characteristic function
$f(t) = \E\,e^{itX}$. We now consider the normalized sum
$$
Z_n = \frac{X_1 + \dots + X_n}{\sqrt{n}},
$$
assuming that $M = M(Z_{n_0})$ is finite. As already discussed 
in Section 2, in this case all $Z_n$ with $n \geq 2n_0$ have 
continuous bounded densities expressed by the inversion formula
$$
p_n(x) = \frac{1}{2\pi} \int_{-\infty}^\infty e^{-itx} f_n(t)\,dt,
\quad x \in \R,
$$
where
$$
f_n(t) = f\Big(\frac{t}{\sqrt{n}}\Big)^n
$$
denotes the characteristic functions of $Z_n$. Equivalently,
\be
p_n(x) = \frac{\sqrt{n}}{2\pi} \int_{-\infty}^\infty
e^{-itx\sqrt{n}}\,f(t)^n\,dt.
\en

Using contour integration, one can cast this formula
in a different form involving the log-Laplace transform 
$K(z) = \log \E\,e^{zX}$ and the saddle point $z_0 = z_0(\tau)$
for the real value $\tau = x/\sqrt{n}$. This is a preliminary
step towards Theorems 1.2 and  1.4. 

As before, let $\E\,e^{\alpha|X|} \leq 2$ with a parameter $\alpha>0$. 

\vskip5mm
{\bf Lemma 6.1} {\sl Let $n \geq 4n_0$. If 
$0 < \ep \leq \frac{\alpha^3}{16}$ and $|\tau| \leq \frac{\ep}{2}$,
then
\be
p_n(x) = \frac{\sqrt{n}}{2\pi} \int_{-\ep}^\ep
\exp\Big\{n (K(z_0 + it) - \tau (z_0 + it))\Big\}\,dt + \theta R_n
\en
with $|\theta| \leq 1$ and
\be
R_n = 5M \exp\Big\{-\frac{n \ep^2}{C n_0 M^2}\Big\}, \quad C = 5200.
\en
}

\vskip2mm
{\bf Proof.} Applying Corollary 3.3, we
get from (6.1) that, for any $\ep \in (0,1]$,
\be
\Big|\,p_n(x) - \frac{\sqrt{n}}{2\pi} \int_{-\ep}^\ep
e^{-itx\sqrt{n}}\,f(t)^n\,dt\Big| \leq 4M
\exp\Big\{-\frac{n\ep^2}{Cn_0 M^2}\Big\}
\en
with $C = 5200$.

Due to the assumption on $\ep$, we may apply Lemma 4.3
which gives
\be
|f(\pm \ep + iy)| \leq e^{-\ep^2/5} \ \ {\rm whenever} 
\ |y| \leq \frac{\ep}{2}.
\en
Assuming for definiteness that $x \geq 0$, we take the rectangle 
contour
$$
L = L_1 + L_2 + L_3 + L_4
$$
with segment parts 

\vskip5mm
$L_1 = [-\ep,\ep]$, 
$L_2 = [\ep,\ep-ih]$, 

$L_3 = [\ep-ih,-\ep-ih]$,
$L_4 = [-\ep-ih,-\ep]$,

\vskip5mm
\noindent
where $h>0$ is chosen to satisfy 
$h \leq \frac{\ep}{2}$. With this choice the complex
numbers $z = t+iy$ with $|t| \leq \ep$, $|y| \leq h$ lie in the
domain of the definition of $K(z)$. Then, by Cauchy's theorem,
\bee
\int_{L_1} e^{-izx\sqrt{n}}\,f(z)^n\,dz + 
\int_{L_2} e^{-izx\sqrt{n}}\,f(z)^n\,dz
 & & \\
 & & \hskip-50mm + \ 
\int_{L_3} e^{-izx\sqrt{n}}\,f(z)^n\,dz + 
\int_{L_4} e^{-izx\sqrt{n}}\,f(z)^n\,dz \, = \, 0.
\ene

Note that in the lower half-plane $z = t - iy$, $0 \leq y \leq h$,
we have $|e^{-izx\sqrt{n}}| = e^{-yx\sqrt{n}} \leq 1$.
Moreover, $|f(z)|$ is bounded away from 1 on $L_2$ and $L_4$
according to (6.5) which gives
$$
\Big| \int_{L_2} \Big| + \Big| \int_{L_4} \Big| \leq 
\ep\,e^{-n\ep^2/5} \leq \frac{1}{16}\,e^{-n\ep^2/5}.
$$
To simplify, note that
$$
4M \exp\Big\{-\frac{n\ep^2}{Cn _0 M^2}\Big\} + 
\frac{1}{16}\,e^{-n\ep^2/5} \leq R_n,
$$
where we used $M \geq \frac{1}{12}$ (cf. Remark 2.3).
Combining this bound with (6.4), we arrive at
\bee
p_n(x)
 & = &
\frac{\sqrt{n}}{2\pi} \int_{L_3}
e^{-izx\sqrt{n}}\,f(z)^n\,dz + \theta R_n \\
 & = &
\frac{\sqrt{n}}{2\pi} \int_{-\ep}^\ep
e^{-i(t - ih)\,x\sqrt{n}}\,f(t-ih)^n\,dt + \theta R_n.
\ene

Using the log-Laplace transform,
let us rewrite the above as a contour integral
$$
p_n(x) = \frac{\sqrt{n}}{2\pi i} \int_{h - i\ep}^{h + i\ep}
\exp\{n (K(z) - \tau z)\}\,dz + \theta R_n
$$
with $\tau = x/\sqrt{n}$ and apply it with 
$h = z_0 = z_0(\tau)$. Due to the requirement
$0 \leq \tau \leq \frac{\ep}{2}$, we have 
$0 \leq \tau \leq \frac{\alpha^3}{32}$ and
$0 \leq z_0 \leq \frac{\alpha^3}{16}$,
according to (5.3), so that Proposition 5.2 is applicable.
After the change of variable, we thus obtain (6.2)-(6.3).
\qed

\vskip5mm
As a next step, let us show that, at the expense of a small error,
the integration in (6.2) may be restricted to the interval
$|t| \leq t_n$ with 
$$
t_n = n^{-1/2}\,\sqrt{8 \log n}, \quad n \geq 4n_0.
$$
This can be achieved under stronger conditions such as
\be
|\tau| \leq \frac{\ep}{80\,M^2 n_0}, \quad 
0 \leq \ep \leq \frac{\alpha^3}{80}.
\en

Indeed, using (5.10) and (5.12) in the representation (6.2), one may
rewrite (6.2) as
$$
p_n(x) = \frac{\sqrt{n}}{2\pi} \int_{-\ep}^\ep
\exp\bigg\{n\, \Big(\sum_{k=2}^\infty \rho_k
\frac{(it)^k}{k!}\Big)\bigg\}\,dt \
e^{n\,(-\frac{1}{2}\,\tau^2 + \tau^3 \lambda(\tau))}
+ \theta R_n,
$$
where $\tau = x\sqrt{n}$ and $\rho_k = K^{(k)}(z_0)$. Equivalently
\be
\frac{p_n(x)}{\varphi(x)} = \frac{\sqrt{n}}{\sqrt{2\pi}} \, 
e^{n\tau^3 \lambda(\tau)} \int_{-\ep}^\ep
\exp\bigg\{n\, \Big(\sum_{k=2}^\infty \rho_k
\frac{(it)^k}{k!}\Big)\bigg\}\,dt + 
\theta R_n e^{x^2/2}.
\en
Here, the new remainder term 
$$
R_n e^{x^2/2} = 
5M \exp\Big\{-\frac{\ep^2 n}{C M^2 n_0} + \tau^2 n\Big\}, \quad
C = 5200,
$$
is still exponentially small with respect to $n$ due to the first
condition in (6.6), which strengthens the assumption 
$|\tau| \leq \frac{\ep}{2}$ in Lemma 6.1 (recall that $M \geq \frac{1}{12}$). 
In this case, the expression in the exponent will be
of order $-\frac{cn\ep^2}{n_0 M^2}$ up to 
an absolute constant $c>0$. Hence, (6.7) yields
\be
\frac{p_n(x)}{\varphi(x)} = \frac{\sqrt{n}}{\sqrt{2\pi}} \, 
e^{n\tau^3 \lambda(\tau)} \int_{-\ep}^\ep
\exp\bigg\{n\, \Big(\sum_{k=2}^\infty \rho_k
\frac{(it)^k}{k!}\Big)\bigg\}\,dt + \theta R_n,
\en
where
\be
R_n = 5M \exp\Big\{-\frac{c\ep^2 n}{M^2 n_0}\Big\}.
\en

Now, by Lemma 4.2,
\be
\frac{1}{2} \leq \rho_2 \leq \frac{3}{2}, \qquad 
|\rho_k| \leq 3 k!\,\Big(\frac{4}{\alpha}\Big)^k \quad (k \geq 3).
\en
It follows that
\bee
{\rm Re}\bigg(\sum_{k=2}^\infty
\rho_k\,\frac{(it)^k}{k!}\bigg)
 & = &
-\rho_2\,\frac{t^2}{2}  + \sum_{k=2}^\infty (-1)^k
\rho_{2k}\,\frac{t^{2k}}{(2k)!} \\
 & \leq &
-\frac{1}{4}\,t^2 + 3 \sum_{k=2}^\infty 
\Big(\frac{4t}{\alpha}\Big)^{2k} \\
 & = &
-\frac{1}{4}\,t^2 + 3\,(16 t^2) \sum_{k=2}^\infty 
\frac{(4t)^{2k-2}}{\alpha^{2k}} 
 \, \leq \, -\frac{1}{8}\,t^2,
\ene
where we used $\alpha<1$ and 
$|t| \leq \ep \leq \frac{1}{80}\,\alpha^3$ so as to bound 
the last sum, according to the second assumption in (6.6). 
Hence, when restricted to $|t| \geq t_n$, the absolute value 
of the integral in (6.8) does not exceed
$$
2 \int_{t_n}^\infty e^{-nt^2/8}\,dt  =
\frac{4}{\sqrt{n}} \int_{\frac{t_n}{2}\sqrt{n}}^\infty e^{-s^2/2}\,ds <
\frac{2\sqrt{2\pi}}{\sqrt{n}}\,e^{-nt_n^2/8} = 
\frac{2\sqrt{2\pi}}{n^{3/2}}.
$$
As a result, assuming the conditions (6.6),
\begin{eqnarray}
\frac{p_n(x)}{\varphi(x)}
 & = &
\frac{\sqrt{n}}{\sqrt{2\pi}} \, e^{n\tau^3 \lambda(\tau)} 
\int_{|t| \leq t_n'} \exp\bigg\{n\, \Big(\sum_{k=2}^\infty
\rho_k \frac{(it)^k}{k!}\Big)\bigg\}\,dt \nonumber \\
 & & + \
2\theta_1 n^{-1} e^{n\tau^3 \lambda(\tau)} + 
\theta_2 R_n,
\end{eqnarray}
where $t_n' = \min(t_n,\ep)$, $|\theta_j| \leq 1$, and where
$R_n$ is now defined in (6.9).

\vskip5mm
\section{{\bf Proof of Theorem 1.4}}
\setcounter{equation}{0}

\vskip2mm
\noindent
As a final step, we need to explore an asymptotic behavior
of the integral in (6.11), where we recall that
$\rho_k = K^{(k)}(z_0)$, $z_0 = z_0(\tau)$ being the saddle
point for $\tau = x/\sqrt{n}$. In view of the conditions
in (6.6), we choose 
$$
\ep = \frac{\alpha^3}{80}, \quad
\tau_0 = \frac{\ep}{80\,M^2 n_0} = \frac{c_0\alpha^3}{M^2\,n_0}
$$
with $c_0 = 1/6400$.
Note that with this choice the definition (6.9) becomes
\be
R_n = 5M \exp\Big\{-\frac{c_1 \alpha^6 n}{M^2 n_0}\Big\},
\en
where $c_1>0$ is an absolute constant.
Suppose that $|\tau| \leq c\tau_0$ with a constant $0 < c \leq 1$
to be chosen later on.

The integrand in (6.11) may be written as
$$
u_n(t) =
\exp\Big\{-n\rho_2\,\frac{t^2}{2} + n\rho_3\,\frac{(it)^3}{6} + 
nv(t)\Big\}
$$
with
$$
v(t) =  \sum_{k=4}^\infty
\rho_k\,\frac{(it)^k}{k!}.
$$

First assume that $n \geq n_1 = \max(4n_0, \ep^{-4})$ which
insures that $t_n' = t_n$ (since $\ep < \frac{1}{80}$). 
As $|t| \leq t_n$, from (6.10) it follows that
\be
nv(t) = O(nt^4) = B\alpha^{-4}\, \frac{(\log n)^2}{n}.
\en
Here and below $B$ denotes a quantity, perhaps different in
different places, bounded by an absolute constant. With this 
convention, since $n \geq  (\frac{80}{\alpha^3})^4$, 
we also have $nv(t) = B$, and by (6.10) with $k=3$,
\be
n\rho_3\,t^3 = B\alpha^{-3}\, \frac{(\log n)^{3/2}}{\sqrt{n}} = B.
\en
So, one may use the Taylor expansion
$e^x = 1 + x + O(x^2)$ in a bounded interval $|x| \leq B$ with
$$
x = n\rho_3\,\frac{(it)^3}{6} + nv(t).
$$
From (7.2)-(7.3) and using again $n \geq B \alpha^{-12}$
together with $\alpha < 1$, we have
$$
(nv(t))^2 = B\alpha^{-8}\, \frac{1}{n} \, \frac{(\log n)^4}{n} = 
\frac{B}{n}
$$
and
$$
(n\rho_3 t^3) \cdot (nv(t)) = B\alpha^{-7}\, 
\frac{1}{n} \, \frac{(\log n)^{7/2}}{\sqrt{n}} = 
\frac{B\alpha^{-2}}{n}.
$$
Since
$$
(n\rho_3 t^3)^2 = B\alpha^{-6}\, \frac{(\log n)^3}{n},
$$
which dominates (7.2) and the previous two expressions, 
we obtain that
\bee
u_n(t) 
 & = &
e^{-n\rho_2 t^2/2 + x} \\
 & = & 
e^{-n\rho_2 t^2/2}\,\big(1 + x + Bx^2\big) \\
 & = &  
e^{-n\rho_2 t^2/2}\,\Big(1 + n\rho_3\,\frac{(it)^3}{6} + 
B\alpha^{-6}\, \frac{(\log n)^3}{n}\Big).
\ene
Hence
$$
\int_{|t| \leq t_n} u_n(t)\,dt = 
\Big(1 + B\alpha^{-6}\, \frac{(\log n)^3}{n}\Big)
\int_{|t| \leq t_n} e^{-n\rho_2 t^2/2}\,\,dt,
$$
and (6.11) is simplified to
\begin{eqnarray}
\frac{p_n(x)}{\varphi(x)}
 & = &
\frac{\sqrt{n}}{\sqrt{2\pi}} \, e^{n\tau^3 \lambda(\tau)} 
\Big(1 + B\alpha^{-6}\, \frac{(\log n)^3}{n}\Big)
\int_{|t| \leq t_n} e^{-n\rho_2 t^2/2}\,\,dt \nonumber \\
 & & + \
2\theta_1 n^{-1} e^{n\tau^3 \lambda(\tau)} + 
\theta_2 R_n.
\end{eqnarray}

Next, one may extend the integration in (7.4) to the whole
real line at the expense of an error not exceeding
\bee
\int_{|t| > t_n} e^{-n\rho_2 t^2/2}\,\,dt
 & = &
\frac{2}{\sqrt{\rho_2 n}} 
\int_{t_n \sqrt{\rho_2 n}}^\infty e^{-s^2/2}\,ds \\
 & < &
\frac{\sqrt{2\pi}}{\sqrt{\rho_2 n}}\,e^{-n\rho_2 t_n^2/2}
 \, \leq \, 
\frac{2\sqrt{\pi}}{\sqrt{n}}\,e^{-n t_n^2/4}
 \, = \, 
\frac{2\sqrt{\pi}}{\sqrt{n}}\,n^{-2},
\ene
where we used $\rho_2 \geq \frac{1}{2}$. The latter bound 
is dominated by $B\alpha^{-6}\, \frac{(\log n)^3}{n}$,
and since the integral over the whole real line is equal to
$\frac{\sqrt{2\pi}}{\sqrt{\rho_2 n}}$, we obtain from (7.4) 
a simpler representation
$$
\frac{p_n(x)}{\varphi(x)} = 
\frac{1}{\sqrt{\rho_2}} \, e^{n\tau^3 \lambda(\tau)} 
\Big(1 + B\alpha^{-6}\, \frac{(\log n)^3}{n}\Big) +
B n^{-1} e^{n\tau^3 \lambda(\tau)} + B R_n.
$$

Here, the first remainder term may be absorbed in the brackets, 
so that this formula is further simplified to
\be
\frac{p_n(x)}{\varphi(x)} = 
\frac{1}{\sqrt{\rho_2}} \, e^{n\tau^3 \lambda(\tau)} 
\Big(1 + B\alpha^{-6}\, \frac{(\log n)^3}{n}  + 
B\, e^{-n\tau^3 \lambda(\tau)} R_n\Big).
\en
Moreover, one may eliminate the factor $e^{-n\tau^3 \lambda(\tau)}$
in front of $R_n$ by choosing a smaller value of $c_1$ in (7.1)
for a proper absolute constant $c$ appearing in the assumption
$|\tau| \leq c\tau_0$. To this end, it will be sufficient to require that
$$
|\psi(\tau)| \leq \frac{c_1 \alpha^6}{2M^2 n_0},
\quad \psi(\tau) = \tau^3 \lambda(\tau).
$$
Recall that, by Proposition 5.5, 
$|\psi(\tau)| \leq 700\,\alpha^{-3}\,|\tau|^3$ in the
interval $|\tau| \leq \frac{1}{64}\,\alpha^3$ (which is larger than
$|\tau| \leq \tau_0$). Hence, the above bound holds true, as long as
\be
|\tau| \leq \frac{c_1^{1/3}}{(1400\,M^2 n_0)^{1/3}}\,\alpha^3.
\en
Since $M \geq \frac{1}{12}$, we have 
$(M^2 n_0)^{1/3} \leq 12^{4/3} M^2 n_0$. Hence (7.6)
may be strengthened to $|\tau| \leq c\tau_0$ with a suitable
constant $c>0$. Under this condition, from (7.5) we thus get
\be
\frac{p_n(x)}{\varphi(x)} = 
\frac{1}{\sqrt{\rho_2}} \, e^{n\tau^3 \lambda(\tau)} 
\Big(1 + B\alpha^{-6}\, \frac{(\log n)^3}{n}  + B R_n\Big),
\en
where $R_n$ is still defined as in (7.1) with a new constant $c_1$.

It is now useful to note that the last error term in this representation
is dominated by the second last one for sufficiently large $n$.
Indeed, using $e^{-y} \leq 2/y^2$ ($y>0$), we have
$$
R_n \leq 2\,\Big(\frac{M^2 n_0}{c_1 \alpha^6 n}\Big)^2 
\leq \alpha^{-6}\, \frac{1}{n},
$$
where the last inequality holds true for 
$n \geq CM^4 n_0^2\, \alpha^{-12}$ with an absolute constant $C>0$.
This condition is slightly stronger than $n \geq n_1$ which was assumed
before. As a result, (7.7) then yields
$$
\frac{p_n(x)}{\varphi(x)} = 
\frac{1}{\sqrt{\rho_2}} \, e^{n\tau^3 \lambda(\tau)} 
\Big(1 + B\alpha^{-6}\, \frac{(\log n)^3}{n}\Big).
$$

It remains to recall Proposition 5.6 according to which
$$
\rho_2^{-1/2} = K''(z_0(\tau))^{-1/2} =  e^{-\mu(\tau)},
$$
and then we arrive at the assertion in Theorem 1.4
(which is a refinement of Theorem 1.2).

Let us also note that the case $2n_0 \leq n < n_1$ is not
interesting, since then $|x| \leq \tau_0 n_1$, and (1.6) holds
true by choosing a suitable constant in $O$ in (1.6).

\vskip5mm
\section{{\bf Proof of Corollary 1.3}}
\setcounter{equation}{0}

\vskip2mm
\noindent
Starting from (5.13) and (5.16), we have
\bee
n\tau^3 \lambda(\tau) - \mu(\tau)
 & = &
\frac{n}{m!}\,\gamma_m \tau^m + O(|\tau|^{m+1}) -
\frac{1}{2(m-2)!}\,\gamma_m \tau^{m-2} + O(|\tau|^{m-1}) \\
 & = &
\frac{\gamma_m}{m!}\,\tau^{m-2}\,\Lambda(\tau).
\ene
Here
$$
\Lambda(\tau) \, = \,
n \tau^2 - \frac{m(m-1)}{2} + n O(\tau^3) + O(\tau) \, \geq \,
\frac{1}{2}\,\Big[n \tau^2 - \frac{m(m-1)}{2}\Big],
$$
which is bounded away from zero, if $|\tau| \leq \tau_1$ for
some constant $\tau_1 > 0$ and $n \tau^2 = x^2 \geq m^2$. 
In this case, (1.6) immediately yields the desired relation (1.7).

In the remaining bounded interval $|x| \leq m$, this argument does not 
work, and it is better to employ the Chebyshev-Edgeworth expansion 
for the correction $\varphi_m(x)$ in (1.1) (which depends on $n$ as well).
In terms of the first non-zero cumulant, (1.1) may be written more
accurately as
$$
p_n(x) = \varphi(x) + 
\frac{\gamma_m}{m!}\,H_m(x)\varphi(x)\, n^{-\frac{m-2}{2}} +
\frac{1}{1 + |x|^m}\,o\big(n^{-\frac{m-2}{2}}\big),
$$
where $H_m(x)$ denotes the Chebyshev-Hermite polynomial of degree $m$.
As a consequence, for any constant $x_0 > 0$,
$$
\sup_{|x| \leq x_0} \frac{|p_n(x) - \varphi(x)|}{\varphi(x)} =
O\big(n^{-\frac{m-2}{2}}\big),
$$
which is stronger than (1.7), since $m$ is even ($m \geq 4$).

\vskip5mm
{\bf Acknowledgment.} We thank the referees
for careful reading and valuable remarks.

\end{document}